\documentstyle[12pt]{amsart}
\makeatletter
\@namedef{pf*}#1{\par%
    \def\proofname{\normalshape#1}\pf\ignorespaces}
\newcommand{\Comp}{\bold C}
\newcommand{\C}{\hat{\bold C}}
\newcommand{\Real}{\bold R}

\newcommand{\Z}{\bold Z}
\newcommand{\bm}[1]{\mbox{\boldmath$#1$}}

\renewcommand{\Re}{\operatorname{Re}}
\renewcommand{\Im}{\operatorname{Im}}
\newcommand{\Res}{\operatornamewithlimits{Res}}

\newcommand{\ord}{\operatorname{Ord}}
\newcommand{\Max}{\operatorname{Max}}

\newcommand{\rank}{\operatorname{rank}}
\newcommand{\Ker}{\operatorname{Ker}}

\newcommand{\Condsum}[2]{\sum\begin{Sb}\scriptscriptstyle #1\\%
                             \scriptscriptstyle #2\end{Sb}}
\newcommand{\Condprod}[2]{\prod\begin{Sb}\scriptscriptstyle #1\\
                             \scriptscriptstyle #2\end{Sb}}

\newcommand{\bg}{\family{cmr}\size{17}{12pt}\selectfont}

\newcommand{\bigzerou}{\smash{\lower1.2ex\hbox{\bg 0}}}

\theoremstyle{plain}
    \newtheorem{Thm}{Theorem}[section]
    \newtheorem{Prop}[Thm]{Proposition}
    \newtheorem{Lemma}[Thm]{Lemma}
    \newtheorem{Cor}[Thm]{Corollary}
    \newtheorem{IntroThm}{Theorem}[section]
        
\theoremstyle{definition}
    
    \newtheorem{Exa}[Thm]{Example}
    \newtheorem{Prob}{Problem}
        
\theoremstyle{remark}
    \newtheorem{Rmk}[Thm]{Remark}
\numberwithin{equation}{section}
\numberwithin{figure}{section}
    \title{An inverse problem of the flux \\ for minimal surfaces}
    \author{Shin Kato}%
    \author{Masaaki Umehara}%
    \author{Kotaro Yamada}%

    \address{(Kato) 
             Department of Mathematics, Faculty of Science, 
             Osaka City University, Osaka 558, JAPAN 
             \newline\indent
             (Umehara) 
             Department of Mathematics, Graduate School of Science, 
             Osaka University, Toyonaka 560, JAPAN 
             \newline\indent
             (Yamada) 
             Department of Mathematics, Faculty of Science,
             Kumamoto University, Kumamoto 860, JAPAN 
            }
    \email{Kato: shinkato@@sci.osaka-cu.ac.jp 
           \newline\indent
           Umehara: umehara@@math.wani.osaka-u.ac.jp 
           \newline\indent
           Yamada: kotaro@@gpo.kumamoto-u.ac.jp}
    \dedicatory{
        Dedicated to Professor Masaru Takeuchi on his sixtieth birthday} 
    \thanks{This research was supported in part by
            Grant-in-Aid for Scientific Research, the Ministry of 
            Education, Sports and Culture, Japan.}
\begin{document}
\maketitle
\begin{abstract}
 For a complete minimal surface in the Euclidean 3-space,
 the so-called flux vector corresponds to each end.
 The flux vectors are balanced, i.e., the sum of those over all ends
 are zero.
 Consider the following inverse problem: For each balanced $n$ vectors,
 find an $n$-end catenoid which attains given vectors as flux.
 Here, an $n$-end catenoid is a complete minimal 
 surface of genus $0$ with ends asymptotic to the catenoids.
 In this paper, the problem is reduced to solving algebraic 
 equation.
 Using this reduction, it is shown that,
 when $n=4$,
 the inverse problem for $4$-end catenoid has 
 solutions for almost all balanced $4$
 vectors.
 Further obstructions for $n$-end catenoids with parallel flux vectors
 are also discussed.
\end{abstract}

\tableofcontents

\section{Introduction}\label{sec:1}

An {\it $n$-end catenoid\/} is 
a complete immersed minimal surface 
of finite total curvature 
which has zero genus and $n$ catenoid ends.
It is considered as a conformal immersion 
    $x\colon{}\C\setminus\{q_1,\ldots,q_n\}\to \Real^3$, 
where 
    $\C:=\Comp\cup\{\infty\}$ 
and $q_1,\ldots,q_n\in\C$.
Jorge-Meeks surfaces \cite{jm} are typical ones.
Recently, 
new examples of $n$-end catenoids have been found by 
\cite{kar}, \cite{l2}, \cite{xu}, \cite{ross1}, \cite{ross2},
\cite{kat} and \cite{uy}.
They contain examples 
with dihedral or Platonic symmetry groups.
We also remark that
for special classes of minimal surfaces 
with catenoid or flat ends,
some systematic approach has been known 
(see \cite{pen}, \cite{xi}, \cite{l1}).

In each end $q_j$ ($j=1,\ldots,n$) of an $n$-end catenoid, 
the {\it flux vector\/} is defined by
\[
    \varphi_j := \int_{\gamma_j} \vec n \, ds, 
\]
where $\gamma_j$ is a curve surrounding $q_j$ 
from the left 
and $\vec n$ the conormal such that 
$(\gamma_j', \vec n)$ is positively oriented.
Each flux vector $\varphi_j$ is proportional 
to the limit normal vector $\nu(q_j)$ 
with respect to the end $q_j$ 
and the scalar $w(q_j):=\varphi_j/4\pi\nu(q_j)$ 
is called the {\it weight\/} of the end $q_j$.
It is well known that 
the flux vectors satisfy 
a \lq\lq balancing" condition 
so called the {\it flux formula} 
\[
    \sum_{j=1}^n\varphi_j = \sum_{j=1}^n 4\pi\,w(q_j)\,\nu(q_j) = 0. 
\]
It should be remarked that $w(q_j)$ may take a negative value.

Therefore, 
we consider an inverse problem
of the flux formula as follows:

\begin{Prob}
    For given unit vectors $\{v_1,\ldots,v_n\}$ in $\Real^3$, 
    and nonzero real numbers $\{a_1,\ldots,a_n\}$ 
    satisfying $\sum_{j=1}^n a_j v_j=0$, 
    is there an $n$-end catenoid 
    $x\colon{}\C\setminus \{q_1,\ldots,q_n\}\to {\Real}^3$ 
    such that $\nu(q_j)=v_j$ 
    and $a_j$ is the weight at the end $q_j$?
\end{Prob}

We remark that 
Kusner also proposed a similar question (see \cite{ross1}).
By the classification of Barbanel~\cite{ba} and Lopez~\cite{l2}, 
we can see that 
the answer for $n\le 3$ is \lq\lq Yes" 
except for the case when two of $\{v_j\}_{j=1}^n$ coincide. 
For $n\ge 4$, 
the first author \cite{kat} gave an explicit formula 
for existence of an $n$-end catenoid 
with prescribed flux 
when $\C\setminus\{q_1,\ldots,q_n\}$ 
is conformally equivalent to 
the image of its Gauss map 
$\bold S ^2\setminus\{v_1,\ldots,v_n\}$.

In this paper, 
we will generalize the formula in \cite{kat},
and show some existence and non-existence results on the problem.
In Section~\ref{sec:2}, we get the following:

\begin{IntroThm}[Theorem~\ref{thm:2-4}]\label{thm:A}
    For any pair $(\bm{v},\bm{a})$ 
    of unit vectors                         
    $\bm{v}=\{v_1,\ldots,\\ v_n\}$ 
    in ${\Real}^3$ 
    and nonzero real numbers $\bm{a}=\{a_1,\ldots,a_n\}$ 
    satisfying $\sum_{j=1}^n a_j v_j=0$, 
    there is an evenly branched $n$-end catenoid 
    $x\colon{}\C\setminus\{q_1,\ldots,q_n\}\to{\Real}^3$ 
    $(q_j\ne\infty)$ 
    such that 
    the induced metric is complete 
    at the end $q_j$,
    $\nu(q_j)=v_j$ 
    and $a_j$ is the weight at the end $q_j$ 
    $(j=1,\ldots,n)$, 
    if and only if 
    there exist complex numbers 
    $b_1,\ldots,b_n$ 
    satisfying the following conditions\rom:
    \[
    \left\{\begin{aligned}
          b_j\displaystyle{\Condsum{k=1}{k\ne j}^nb_k\frac{p_k-p_j}{q_k-q_j}}
            &=a_j \\
          b_j\displaystyle{
                \Condsum{k=1}{k\ne j}^n b_k\frac{\bar p_jp_k+1}{q_k-q_j}}
            &=0 
    \end{aligned}\right.
    \qquad\qquad 
    (j=1,\ldots,n),
    \]
    where $p_j:=\sigma(v_j)$, 
    $\sigma:\bold S ^2\to\C$ is the stereographic projection, 
    and we assume $p_j\ne\infty$.

    Moreover, 
    the immersion $x$ has no branch points 
    if and only if the resultant $\Psi(P(z),Q(z))$ 
    of the polynomials $P(z)$ and $Q(z)$ 
    $($defined by \eqref{Def:P} and \eqref{Def:Q}$)$ 
    does not vanish.
\end{IntroThm}

We note here that 
the flux formula holds 
even if the surface allows branch points 
(see Remark~\ref{rmk:2-8}).

In the case when an $n$-end catenoid 
has the same symmetry as its flux data, 
the construction is reduced to a routine work 
by virtue of our theorem, 
and one can construct all of the known examples 
(cf. \cite{kar}, \cite{l2}, \cite{xu}, 
\cite{ross1}, \cite{ross2}, \cite{kat}, \cite{uy}, etc.) 
and far more new examples (cf. \cite{kuy3}). 

We also remark here that 
an $n$-end catenoid does not always have 
the symmetry of its flux data. 
In fact, there exists a flux data $(\bm{v},\bm{a})$ 
such that any corresponding surface 
does not have the same symmetry 
as $(\bm{v},\bm{a})$ 
(see Example \ref{exa:3-7}(iii)). 

On the other hand, 
for a certain flux data, 
there are no $n$-end catenoids realizing it. 
Indeed, there are no $n$-end catenoids with the flux data 
$(\bm{v},\bm{a})$ satisfying one of the following conditions: 
\begin{enumerate}
  \item $v_1=v_2=v_3=\cdots=v_n$; 
  \item $-v_1=-v_2=v_3=\cdots=v_n$; 
  \item $-v_1=v_2=\cdots=v_n$ 
        and $\sum_{j=2}^{n-1}\sum_{k=j+1}^na_ja_k\ne 0$; 
  \item $(n=4)$ $-v_1=v_2$ and $v_3=v_4\ne \pm v_1$. 
\end{enumerate}

The first condition is well-known, 
and the third condition follows from the genus zero case 
of the second compatibility condition in \cite{per}. 
The fourth condition is new. 
These four obstructions are easily obtained 
as a corollary of Theorem~A. 

It is interesting to observe that, 
when the equality 
$\sum_{j=2}^{n-1}\sum_{k=j+1}^na_ja_k=0$ 
holds in (3) above, 
$n$-end catenoids which allow the deformation 
described by Lopez-Ros~\cite{lr} 
can be constructed 
(see Examples~\ref{exa:4-8},\ref{exa:4-9} and \ref{exa:4-10}).

In spite of the above non-existence results, 
it seems that generic flux data are free of additional obstructions.
We demonstrate it for $n=4$, 
and get the following:

\begin{IntroThm}[Theorems~\ref{thm:3-3} and \ref{thm:3-6}]\label{thm:C}
    For almost all pair $(\bm{v},\bm{a})$ of 
    unit vectors $\bm{v}=\{v_1,v_2,v_3,v_4\}$ in ${\Real}^3$ 
    and nonzero real numbers $\bm{a}=\{a_1,a_2,a_3,a_4\}$ 
    satisfying $\sum_{j=1}^4a_jv_j =0$, 
    there is a \rom(non-branched\rom) $4$-end catenoid 
    $x\colon{}\C\setminus\{q_1,q_2,q_3,q_4\}\to{\Real}^3$ 
    such that $\nu(q_j)=v_j$ 
    and $a_j$ is the weight at the end $q_j$ 
    $(j=1,2,3,4)$, 
    where $\nu$ is the Gauss map of $x$.
    Moreover, 
    the number of such $x$ is at most $4$ 
    up to rigid motions in ${\Real}^3$.
    In particular, 
    there exist $4$-end catenoids 
    with no symmetric properties.
\end{IntroThm}

We remark here that $n=4$ is the smallest number such that 
$n$-end catenoids have various conformal structures 
and that there exist 
mutually non-congruent $n$-end catenoids 
with the same flux data. 
Indeed, the upper estimate in the theorem above is sharp 
(see Example~\ref{exa:3-7} and Figure~\ref{fig:4noid2}). 

To prove the first part of Theorem~\ref{thm:C}, 
we give an explicit algorithm to construct $4$-end catenoids 
with the prescribed flux 
by reducing it to solve a certain algebraic equation 
of degree $4$. 
However, 
to treat the case when $n\ge5$, 
we shall have to do more complicated analysis 
(cf. \cite{kuy2}). 

\vspace{0.5\baselineskip}

The authors are very grateful to Dr. Wayne Rossman 
for valuable discussions and encouragements. 
They also thank to Prof. Osamu Kobayashi and Dr. Shin Nayatani 
for useful comments. 

\renewcommand{\labelenumi}{\mbox{\rm (\roman{enumi})}}
\renewcommand{\theenumi}{\mbox{\rm (\roman{enumi})}}


\section{Reduction of the problem}\label{sec:2}

For a positive integer $n$, 
we fix a Riemann surface
\[
     M^2=\C\setminus\{q_1,\ldots,q_n\},
\]
where $\C:=\Comp\cup\{\infty\}$ 
and $q_1,\ldots,q_n$ are mutually distinct points. 
A pair $(g,\omega)$ 
of a meromorphic function $g$ 
and non-vanishing meromorphic $1$-form $\omega$ on $\C$ 
is called {\it Weierstrass data}.
By the Enneper-Weierstrass representation formula,
the map defined by
\begin{equation}\label{eq:2-1}
     x:=\Re
        \left(
            \int_{z_0}^z(1-g^2)\omega,
            \int_{z_0}^z i(1+g^2)\omega,
            \int_{z_0}^z 2g\omega
        \right) 
\end{equation}
is a branched minimal immersion 
on the universal covering on $M^2$,
where $z_0$ is a fixed point on $M^2$.
Any isolated degenerate point of 
the induced metric
\begin{equation}\label{ds2}
    ds^2=(1+|g|^2)^2|\omega|^2
\end{equation}
is corresponding to a branch point
of the map. 
The branched minimal immersion $x$ 
is single-valued on $M^2$ 
if and only if
\begin{equation}\label{Res}
    \left\{
    \begin{aligned}
        \Res_{z=q_j}(g\omega)&\in\Real,\\
        \Res_{z=q_j}(\omega)&=
            -\overline{\Res_{z=q_j}(g^2\omega)},
    \end{aligned}
    \right.
    \qquad\qquad
    (j=1,\ldots,n),
\end{equation}
where ${\Res}_{z=q_j}$ is the residue at $z=q_j$.
Moreover, 
if $ds^2$ is a complete Riemannian metric on $M^2$,
then the map $x$ is a complete minimal immersion on $M^2$ 
with finite total curvature.
Conversely,
any complete conformal minimal immersion 
$x\colon{}M^2\to{\Real}^3$ 
with finite total curvature 
is constructed in such manner 
from the following Weierstrass data: 
\begin{align}
     g      &= \frac{\partial x_3}{\partial x_1-i\partial x_2}, 
               \label{Wei1} \\
     \omega & =\partial x_1-i\partial x_2,              
               \label{Wei2}
\end{align}
where the function $g$ is 
the stereographic projection of the Gauss map.

In this section, 
we rewrite the condition \eqref{Res}
into purely algebraic ones.
We remark that 
the second fundamental form 
of the minimal immersion $x$ 
is expressed by 
    $-\omega\cdot dg-\overline{\omega\cdot dg}$ 
and its $(2,0)$-part $\omega\cdot dg$ 
is called the {\it Hopf differential\/} 
of the immersion $x$. 
The end $q_j$ is called a {\it catenoid end\/} 
if the end is asymptotic to a catenoid 
by a suitable homothety, 
that is, 
the Gauss map has no branch point at the each end $q_j$ 
and the Hopf differential $\omega\cdot dg$ 
has a pole of order $-2$ at $q_j$.
The minimal immersion $x$ is called 
an {\it $n$-end catenoid\/}
if all ends $q_1,\ldots,q_n$ are catenoid ends.
We also use a terminology {\it branched $n$-end catenoid\/}
when the induced metric allows 
at most finite degenerate points.
In particular, 
we call a branched $n$-end catenoid is {\it evenly branched\/} 
if all of its branch points are of even order.
We regard non-branched $n$-end catenoids as 
special cases of evenly branched $n$-end catenoids.

First we prepare the following lemma:

\begin{Lemma}\label{lem:2-1}
    Let $x\colon{}M^2\to {\Real}^3$ be a branched $n$-end catenoid.
    Then the degree of its Gauss map  is at most $n-1$ and
    by a suitable motion in ${\Real}^3$, the Weierstrass data
    given by \eqref{Wei1} and \eqref{Wei2} are taken 
    to be satisfying the following conditions\rom :
    \begin{enumerate}
    \item  $\omega$ has poles of order $-2$ on $\{q_1,\ldots,q_n\}$.
    \item  $g$ has no poles and branch points on $\{q_1,\ldots,q_n\}$.
    \end{enumerate} 
    Moreover, $x$ has no branch points if and only if the
    degree of the Gauss map is $n-1$.
\end{Lemma}
\begin{pf}
    By a suitable motion in ${\Real}^3$, 
    we may assume that 
    $g$ has no poles on $\{q_1,\ldots,q_n\}$.
    We apply the relation 
    \begin{equation}\label{index}
     \sum_{z\in Z(\omega)}\ord_z(\omega)
         +\sum_{z\in S(\omega)}\ord_z(\omega)
          = -\chi(\C)=-2, 
    \end{equation}
    where $Z(\omega)$ and $S(\omega)$ are 
    the set of zeros and the set of poles respectively.
    The assumption of catenoid ends yields that
    the Hopf differential $\omega\cdot dg$ 
    has a pole of order $-2$ 
    and $dg$ has no zero at each end $q_j$.
    So $\omega$ has exactly order $-2$ at each end $q_j$.
    Therefore we have that
    \begin{equation}\label{poles}
         \sum_{z\in S(\omega)}\ord_z(\omega)=-2n.
    \end{equation}
    By \eqref{index} and \eqref{poles}, 
    we have
    \begin{equation}\label{zeros}
         \sum_{z\in Z(\omega)}\ord_z(\omega)=2n-2.
    \end{equation}

    On the other hand, 
    since $g$ has no pole at each end, 
    any pole of $g$ must be a zero of $\omega$ by \eqref{ds2}.
    In particular, 
    the inequality
    \begin{equation}\label{metric}
         \ord_z(\omega)\ge-2\ord_z(g)
         \qquad 
         (z\in Z(\omega)).
    \end{equation}
    holds.
    Since the degree $\deg(g)$ of Gauss map is given by
    \[
         \deg(g)=-\sum_{z\in S(g)}\ord_z(g), 
    \]
    we have
    \begin{equation}\label{deg}
         \sum_{z\in Z(\omega)}
         \ord_z(\omega)\ge 2\deg(g).
    \end{equation}
    By \eqref{zeros} and \eqref{deg}, 
    we get
    \begin{equation}\label{degG}
         \deg(g)\le n-1. 
    \end{equation}
    Here,
    $x$ has no branch points 
    if and only if \eqref{metric} is an equality, 
    and hence the equality of \eqref{degG} holds 
    if and only if $x$ has no branch points.
\end{pf}

\begin{Lemma}\label{lem:2-2}
    Let $M^2=\C\setminus\{q_1,\ldots,q_n\}$.
    Let $g$ be a meromorphic function 
    and $\omega$ a meromorphic $1$-form on $\C$ 
    satisfying the conditions \rom{(i)} and \rom{(ii)} 
    of the Lemma~\ref{lem:2-1}.
    Set $p_j:=g(q_j)$.
    Assume $q_j\ne\infty$ and $p_j\ne\infty$ $(j=1,\ldots,n)$.
    Then the symmetric tensor
    \begin{equation}\tag{2.2}
         ds^2=(1+|g|^2)^2|\omega|^2 
    \end{equation}
    is a complete Riemannian metric on $M^2$
    if and only if
    there exist two polynomials
    \begin{align}
        Q(z) & =  
            \sum_{j=1}^nb_j\Condprod{k=1}{k\ne j}^n(z-q_k), 
        \label{Def:Q} \\
        P(z) & = 
            \sum_{j=1}^np_jb_j\Condprod{k=1}{k\ne j}^n(z-q_k) 
        \label{Def:P}
    \end{align}
    $(b_1,\ldots,b_n\in{\Comp})$ 
    satisfying the following properties\rom:
    \begin{enumerate}
    \item    $\Max\{\deg(P),\deg(Q)\}=n-1$.
    \item    $P(z)$ and $Q(z)$ are irreducible.
    \item    $g(z)=P(z)/Q(z)$ and 
             $\omega(z)=-\{\sum_{j=1}^nb_j/(z-q_j)\}^2 dz$.
    \end{enumerate} 
\end{Lemma}

\begin{pf}
    We suppose that $ds^2$ is a complete Riemannian metric on $M^2$.
    Since $ds^2$ is positive definite on $M^2$,
    we have from \eqref{ds2} that 
    the inequality \eqref{metric} turns to be an equality
    \begin{equation}\label{mmetric}
         \ord_z(\omega)=-2\ord_z(g) 
         \qquad  
         (z \in Z(\omega)\,).
    \end{equation}
    By the same argument in the proof of the previous lemma, 
    we have
    \begin{equation}\label{degGG}
         \deg(g)=n-1.
    \end{equation}
    Since $\omega$ has only poles of order $-2$, 
    poles and zeros of $\omega$ are all even order.
    Thus $\sqrt{\omega}$ is defined 
    as a meromorphic section 
    of the half canonical line bundle 
    and has poles of order $-1$ on $P(\omega)$.
    Since $\omega$ has no pole at infinity, 
    there exist complex numbers 
    $b_1,\ldots,b_n\in{\Comp}$ such that
    \[
         \sqrt{\omega}
              = i\left(
                    \sum_{j=1}^n\frac{b_j}{z-q_j}
                 \right)
                 \sqrt{dz}.
    \]

    Now we set
    \begin{align}
         R(z)&:=\prod_{j=1}^n(z-q_j), \label{Def:R}\\
         R_j(z)&:=\Condprod{k=1}{k\ne j}^n(z-q_k). \label{Def:Rj}
    \end{align}
    Then we have $\omega=-\left(\sum_{j=1}^n b_jR_j(z)\right)^2/R(z)^2$.
    Hence, 
    by \eqref{mmetric}, 
    $g$ can be written as
    \[
         g=\frac{P(z)}{\sum_{j=1}^n b_jR_j(z)},
    \]
    where $P(z)$ is a polynomial of order $n-1$.
    Clearly, 
    $Q(z)$ defined by \eqref{Def:Q} 
    satisfies $Q(z)=\sum_{j=1}^n b_jR_j(z)$.
    Moreover, 
    we have \eqref{Def:P} since $g(q_j)=p_j$.
    By \eqref{degGG}, 
    $P(z)$ and $Q(z)$ are irreducible 
    and $\Max\{\deg(P),\deg(Q)\}=n-1$. 
    On the other hand, 
    the symmetric tensor $ds^2$ 
    induced from such two polynomials $P(z)$ and $Q(z)$ 
    by \eqref{ds2}
    is obviously a complete Riemannian metric on $M^2$.
\end{pf}

The following proposition reduces 
the conditions \eqref{Def:Q} and \eqref{Def:P}
to a purely algebraic condition, 
which plays essential roles in this paper: 

\begin{Prop}\label{prop:2-3}
    Let $n\ge 2$ be an integer and
    $q_1,\ldots,q_n,p_1,\ldots,p_n,b_1,\ldots,b_n$ complex numbers.
    Set $g(z):=P(z)/Q(z)$ 
    and $\omega(z):=-\{\sum_{j=1}^n b_j/(z-q_j)\}^2 dz$,
    where $P(z)$ and $Q(z)$ are polynomials 
    defined by \eqref{Def:Q} and \eqref{Def:P} respectively.
    Then the branched minimal immersion 
    \[
         x = \Re\left(
                \int_{z_0}^z(1-g^2)\omega,
                \int_{z_0}^z i(1+g^2)\omega,
                \int_{z_0}^z 2g\omega\right)
    \]
    is single-valued on the Riemann surface 
    $M^2=\C\setminus\{q_1,\ldots,q_n\}$ 
    if and only if the following conditions hold\rom:
    \begin{equation}\label{Cond}
        \left\{
        \begin{aligned}
          b_j\displaystyle{\Condsum{k=1}{k\ne j}^nb_k\frac{p_j-p_k}{q_j-q_k}}
           &\in {\Real} \\
          b_j
            \displaystyle{
                \Condsum{k=1}{k\ne j}^nb_k\frac{\bar p_jp_k+1}{q_j-q_k}}
           &=  0
        \end{aligned}
        \right.
        \qquad\qquad 
        (j=1,\ldots,n).
    \end{equation}
\end{Prop}

\begin{pf}
    By \eqref{Def:Q} and \eqref{Def:P}, 
    one can easily get the following identities
    \begin{align*}
        \Res_{z=q_j}(\omega) & = 
            -2b_j\Condsum{k=1}{k\ne j}^n\frac{b_k}{q_j-q_k},\\
        \Res_{z=q_j}(g\omega) & = 
            -b_j\Condsum{k=1}{k\ne j}^n b_k\frac{p_j+p_k}{q_j-q_k},\\
         \Res_{z=q_j}(g^2\omega) & = 
            -2b_jp_j\Condsum{k=1}{k\ne j}^n b_k\frac{p_k}{q_j-q_k}.
    \end{align*}
    Thus the conditions \eqref{Res} can be rewritten as
{\allowdisplaybreaks
    \begin{alignat}{2}
        &\left\{\begin{aligned}
          b_j\displaystyle{\Condsum{k=1}{k\ne j}^nb_k\frac{p_j+p_k}{q_j-q_k}}
            &\in{\Real}\\
          b_j\displaystyle{\Condsum{k=1}{k\ne j}^n\frac{b_k}{q_j-q_k}}
            &= -\displaystyle{
                \overline{b_jp_j\Condsum{k=1}{k\ne j}b_k\frac{p_k}{q_j-q_k}}}
        \end{aligned}\right.
        \qquad\qquad
        &&(j=1,\ldots,n).
        \tag{2.3$'$}
        \\
    \intertext{If we set}
        \label{AB}
        &\left\{\begin{aligned}
          A_j &:= 
            b_j\displaystyle{\Condsum{k=1}{k\ne j}^n\frac{b_k}{q_j-q_k}}\\
          B_j &:= 
            b_j\displaystyle{\Condsum{k=1}{k\ne j}^nb_k\frac{p_k}{q_j-q_k}}
         \end{aligned}\right.
        \qquad\qquad 
        &&(j=1,\ldots,n),\\
    \intertext{
        then \eqref{Res} is equivalent 
        to the following condition:}
        \label{C}
        &\left\{\begin{aligned}
              p_jA_j+B_j\in {\Real} \\
              A_j+\bar p_j\overline B_j = 0
        \end{aligned}\right.          
        \qquad\qquad
        &&(j=1,\ldots,n).\\
    \intertext{
        It can be easily seen that 
        \eqref{C} implies that 
        $p_jA_j$ and $B_j$ are both real numbers.
        Hence \eqref{C} reduces 
        to the following condition:}
        \label{CC}
        &\left\{\begin{aligned}
              p_jA_j-B_j\in {\Real} \\
              A_j+\bar p_jB_j = 0 
        \end{aligned}\right.
        \qquad\qquad
        &&(j=1,\ldots,n),
    \end{alignat}
}                                   
    which is equivalent to 
    the desired condition \eqref{Cond}.
\end{pf}

Next we consider the flux formula
on $n$-end catenoid.
We fix a Riemann surface
\[
    M^2=\C\setminus\{q_1,\ldots,q_n\}, 
\]
where $q_1,\ldots,q_n$ are mutually distinct points. 
Let $x\colon{}M^2\to{\Real^3}$ 
be a branched $n$-end catenoid.
Then the {\it flux vector\/} at an end $q_j$ is defined by
\begin{equation}\label{DefF}
     \varphi_j := \int_{\gamma_j} \vec n\, ds, 
\end{equation}
where $\gamma_j$ is a circle 
surrounding $q_j$ from the left, 
and $\vec n$ the conormal such that 
$(\gamma_j',\vec n)$ is positively oriented.
The flux vector is independent 
of  choice of a circle $\gamma_j$.
Each flux vector $\varphi_j$ is 
proportional to the limit normal vector $\nu(q_j)$ 
with respect to the end~$q_j$, 
and the real number $a_j:=\varphi_j/4\pi\nu(q_j)$ 
is called the {\it weight\/} of the end $q_j$.
One can easily verify that 
the Hopf differential $\omega\cdot dg$ 
has the following Laurent expansion 
at each end $q_j$:
\begin{equation}\label{Hopf}
    \omega\cdot dg=\left\{\frac{a_j}{(z-q_j)^2}+\cdots\right\}dz^2,
\end{equation}
where $a_j$ is the weight at the end $q_j$.
By Lemma~\ref{lem:2-1}, 
we may assume that 
$g$ has no poles 
and $\omega$ has poles of order $-2$ on $\{q_1,\ldots,q_n\}$.
Then, 
by Lemma~\ref{lem:2-2}, 
there exist complex numbers 
$p_1,\ldots,p_n$ and $b_1,\ldots,b_n$ 
such that 
\newline
$g(z)=P(z)/Q(z)$ and $\omega(z)=-\{\sum_{j=1}^n b_j/(z-q_j)\}^2 dz$,
where $P(z)$ and
$Q(z)$ are polynomials 
defined by \eqref{Def:Q} and \eqref{Def:P} respectively.
Then, 
by \eqref{Hopf}, 
we have the identities
\begin{equation}\label{a}
     a_j=b_j\Condsum{k=1}{k\ne j}^nb_k\frac{p_j-p_k}{q_j-q_k} 
    \qquad (j=1,\ldots,n).
\end{equation}
We remark that the reality of $a_j$ 
follows from the conditions \eqref{Cond} and \eqref{a}.
Since the limit normal vector $\nu(q_j)$ 
with respect to the end $q_j$ 
is expressed by
\[
 \nu(q_j) = \left(\frac{2\Re(p_j)}{|p_j|^2+1},
                  \frac{2\Im(p_j)}{|p_j|^2+1},
                  \frac{|p_j|^2-1}{|p_j|^2+1}
            \right), 
\]
as the inverse stereographic image of $p_j$,
the flux formula stated 
in the introduction is rewritten as
\begin{equation}\label{Flux}
    \sum_{j=1}^na_j\frac{|p_j|^2-1}{|p_j|^2+1}=0, \qquad
    \sum_{j=1}^na_j\frac{\bar p_j}{|p_j|^2+1}=0. 
\end{equation}

As an application of Lemmas~\ref{lem:2-1}, 
\ref{lem:2-2} and Proposition~\ref{prop:2-3}, 
we have the following reduction theorem 
for the inverse problem of the flux formula:

\begin{Thm}\label{thm:2-4}
    Let $n\ge 2$ be an integer,
    $q_1,\ldots,q_n,p_1,\ldots,p_n,b_1,\ldots,b_n$ complex numbers, 
    and $a_1,\ldots,a_n$ real numbers.
    Set $g(z):=P(z)/Q(z)$ 
    and 
    $\omega(z):=-\{\sum_{j=1}^n b_j/(z-q_j)\}^2 dz$,
    where $P(z)$ and $Q(z)$ are polynomials 
    defined by \eqref{Def:Q} and \eqref{Def:P} respectively. 
    Then the map 
    \[
         x:= \Re\left(
                \int_{z_0}^z(1-g^2)\omega,
                \int_{z_0}^z i(1+g^2)\omega,
                \int_{z_0}^z 2g\omega
              \right)  
    \]
    is an evenly branched $n$-end catenoid defined on  
    $M^2=\C\setminus\{q_1,\ldots,q_n\}$ 
    and has the flux vector
    \[
         \varphi_j = 4\pi a_j\cdot\left(\frac{2\Re (p_j)}{|p_j|^2+1},
                           \frac{2\Im (p_j)}{|p_j|^2+1},
                           \frac{|p_j|^2-1}{|p_j|^2+1}
                     \right)
    \]
    at each end $q_j$ 
    if and only if the following condition holds\rom:
    \begin{equation}\label{Reduction2}
    \left\{\begin{aligned}
          b_j\displaystyle{\Condsum{k=1}{k\ne j}^nb_k\frac{p_k-p_j}{q_k-q_j}}
            &=a_j \\
          b_j\displaystyle{
                \Condsum{k=1}{k\ne j}^n b_k\frac{\bar p_jp_k+1}{q_k-q_j}}
            &=0 
    \end{aligned}\right.
    \qquad\qquad 
    (j=1,\ldots,n).
    \end{equation}
    Moreover, 
    suppose $\Max\{\deg(P),\deg(Q)\}=n-1$, 
    and $P(z)$ and $Q(z)$ are irreducible.
    Then $x$ has no branch points 
    and is an $n$-end catenoid.
    Conversely, 
    any evenly branched $n$-end catenoid 
    is constructed in such manner. 
\end{Thm}

\begin{pf}
    We set real numbers $A_j$ and $B_j$ by \eqref{AB}.
    Then, 
    \eqref{Reduction2} is rewritten as
    \[
        \left\{\begin{aligned}
          p_jA_j-B_j &= a_j \\
          -(A_j+\bar p_jB_j) &= 0 
        \end{aligned}\right.
        \qquad\qquad
        (j=1,\ldots,n).
    \]
    Hence the first assertion follows immediately 
    from \eqref{CC} and \eqref{a}.
    Conversely, we fix an evenly branched $n$-end catenoid with 
    Weierstrass data $(g,\omega)$. Then the order of $\omega$ is
    even everywhere. Thus $\sqrt{\omega}$ is defined as a meromorphic section
    of the half canonical bundle. By the same argument of the proof
    of Lemma~\ref{lem:2-2}, we have the expression 
    $g(z):=P(z)/Q(z)$ 
    and $\omega(z):=-\{\sum_{j=1}^n b_j/(z-q_j)\}^2 dz$,
    where $P(z)$ and $Q(z)$ are polynomials 
    defined by \eqref{Def:Q} and \eqref{Def:P} respectively. 
    This proves the second assertion.
\end{pf}

\begin{Rmk}
    The Weierstrass data of 
    any evenly branched $n$-end catenoid 
    have the form as in Proposition~\ref{prop:2-3}.
    Therefore, 
    it is branched if and only if 
    the resultant 
    of $P(z)$ and $Q(z)$ 
    does not vanish.
    This algebraic equation is 
    expected to have zeros of codimension $1$.
    Indeed it is true in the case $n=4$ 
    as we will see in Section~\ref{sec:3}.
\end{Rmk}

\begin{Rmk}\label{rmk:2-9}
When $q_j=rp_j$\,\,$(j=1,\ldots,n)$, 
Theorem 2.4 reduces to the results in the first author \cite{kat}. 
In this case, 
the system \eqref{Reduction2} reduces to 
\[
    \left\{
    \begin{array}{l}
        \displaystyle{\frac{1}{r} b_j 
                      \Condsum{k=1}{k\ne j}^n b_k} 
         = a_j \\ 
        \displaystyle{\frac{1}{r} b_j 
                      \Condsum{k=1}{k\ne j}^n 
                       b_k \frac{\overline{p_j}p_k+1}{p_k-p_j}} 
         = 0 
    \end{array}
    \right.
    \qquad\qquad 
    (j=1,\ldots,n). 
\]
Moreover, 
the surface has no branch points 
if and only if $\sum_{j=1}^n b^j\ne 0$. 
Many known and new examples of $n$-end catenoids 
can be constructed from this formula, 
and also from our formula \eqref{Reduction2} 
(cf. \cite{kuy3}). 
\end{Rmk}

\begin{Rmk}
    When $p_n\ne p_j$ $(j=1,\ldots,n-1)$, 
    we may assume $p_n=q_n=\infty$ without loss of generality.
    Under this assumption, 
    since $p_j\ne\infty$ $(j=1,\ldots,n-1)$ holds automatically, 
    it is easy to see that 
    the equation \eqref{Reduction2} can be rewritten 
    as the following version:
    \begin{equation}\label{Eq}
    \left\{
    \begin{aligned}
     &    b_j\left(
            \displaystyle{
            \Condsum{k=1}{k\ne j}^{n-1}b_k
            \frac{p_k-p_j}{q_k-q_j}+b_n}
            \right)=a_j 
            \qquad (j=1,\ldots,n-1), \\
     &    b_n\sum_{k=1}^{n-1}b_k=a_n, \\
     &    b_j\left(
            \displaystyle{
            \Condsum{k=1}{k\ne j}^{n-1}b_k
            \frac{\bar p_jp_k+1}{q_k-q_j}+\bar p_jb_n}
            \right)=0 
            \qquad (j=1,\ldots,n-1), \\
     &    -b_n\sum_{k=1}^{n-1}p_kb_k=0. 
    \end{aligned}
    \right.
    \end{equation}

    In this situation, 
    the polynomials $P(z)$, $Q(z)$ and $R(z)$ 
    are replaced naturally as follows:
    \begin{align*}
     P(z) & = 
      \sum_{j=1}^{n-1}p_jb_j\Condprod{k=1}{k\ne j}^{n-1}(z-q_k)
      -b_n\prod_{k=1}^{n-1}(z-q_k), \\
     Q(z) & =  
      \sum_{j=1}^{n-1}b_j\Condprod{k=1}{k\ne j}^{n-1}(z-q_k), \\
     R(z) & =  \prod_{k=1}^{n-1}(z-q_k).
    \end{align*}
\end{Rmk}

By easy calculation, 
we get the following

\begin{Cor}\label{cor:2-7}
    The assertion of Theorem~\ref{thm:2-4} holds
    even if we replace the condition \eqref{Reduction2} 
    by the following condition:
    \begin{equation}\label{Red}
    \left\{\begin{aligned}    
          b_j\displaystyle{\Condsum{k=1}{k\ne j}^nb_k\frac{1}{q_j-q_k}}
            &= 
                \displaystyle{a_j\frac{\bar p_j}{|p_j|^2+1}}\\
          b_j\displaystyle{\Condsum{k=1}{k\ne j}^nb_k\frac{p_j+p_k}{q_j-q_k}}
            &= 
                \displaystyle{a_j\frac{|p_j|^2-1}{|p_j|^2+1}} 
    \end{aligned}\right.
        \qquad\quad(j=1,\ldots,n).
    \end{equation}
\end{Cor}

\begin{Rmk}\label{rmk:2-8}
    Summing up the equations \eqref{Red}
    for $j=1,\ldots,n$, 
    we get the flux formula \eqref{Flux} 
    for any evenly branched $n$-end catenoids. 
    However, 
    the flux formula itself is still valid for minimal syrfaces 
    with branch points of odd order. 
    In fact, by straightforward calculation, 
    the flux vector defined by \eqref{DefF} is written as 
    \[
         \varphi_j = -\Im\left(
                \int_{\gamma_j} (1-g^2)\omega,
                \int_{\gamma_j} i(1+g^2)\omega,
                \int_{\gamma_j} 2g\omega
              \right) . 
    \]
    The flux formula is obvious from the point of view. 
\end{Rmk}

\renewcommand{\labelenumi}{\mbox{\rm (\arabic{enumi})}}
\renewcommand{\theenumi}{\mbox{\rm (\arabic{enumi})}}


\section{4-end catenoids of generic type}\label{sec:3}

Let 
    $x\colon{}\C\setminus\{q_1,\ldots,q_n\}\to {\Real}^3$ 
be an {\it $n$-end catenoid}, 
and set $v_j:=\nu(q_j)$ ($j=1,\ldots,n$), 
where $\nu$ is the Gauss map.
Then, 
as we saw in the previous sections, 
a family $\bm{v}=\{v_1,\ldots,v_n\}$ 
of $n$ unit vectors in ${\Real}^3$ 
must satisfy the following condition.
\begin{equation}\tag{F.$n$}\label{fn}
     \sum_{j=1}^na_jv_j=0 
     \qquad
     \text{for some nonzero real numbers}
     \quad
     a_1,\ldots,a_n.
\end{equation}
Now, 
we classify arrangements of 
$\bm{v}=\{v_1,\ldots,v_n\}$ 
to the following three types:
\begin{description}
\item[{\makebox[4.5em][l]{TYPE I}}]
    $\bm{v}=\{v_1,\ldots,v_n\}$ satisfies \eqref{fn} and 
    dim$\langle v_1,\ldots,v_n\rangle=1$.
\item[{\makebox[4.5em][l]{TYPE II}}]
    $\bm{v}=\{v_1,\ldots,v_n\}$ satisfies \eqref{fn} and 
    dim$\langle v_1,\ldots,v_n\rangle=2$.
\item[{\makebox[4.5em][l]{TYPE III}}]
    $\bm{v}=\{v_1,\ldots,v_n\}$ satisfies \eqref{fn} and 
    dim$\langle v_1,\ldots,v_n\rangle=3$.
\end{description}

We call a (branched) $n$-end catenoid is of TYPE I (resp.~II, III), 
if $\bm{v}$ is of TYPE I (resp.~II, III).

The following facts are already known 
(e.g.~\cite{m}, \cite{ba}, \cite{l2}, \cite{kat}): 
\begin{enumerate}\it
\item   There are no $1$-end catenoids.
\item   Any $2$-end catenoid is the catenoid. 
            Of course, 
        it is of TYPE I\@.
\item   There are no $3$-end catenoids of TYPE I\@.
        Consequently, 
        any $3$-end catenoid is of TYPE II\@. 
        More precisely, 
        for any unit vectors $v_1,v_2,v_3$ of TYPE II and 
        nonzero real numbers $a_1,a_2,a_3$ 
        satisfying $\sum_{j=1}^3a_jv_j=0$, 
        there exists an essentially unique $3$-end catenoid 
        $x\colon{}\C\setminus\{q_1,q_2,q_3\}\to {\Real}^3$ 
        which satisfies 
        $\nu(q_j)=v_j$ and the weight $w(q_j)=a_j$ 
        $(j=1,2,3)$.
\end{enumerate}
From these results, 
the moduli of at most $3$-end catenoids 
is understood completely.

In this section, 
we restrict our attention to 
$4$-end catenoids of TYPE III 
that is a generic type.
In particular, 
we give some upper estimates 
for the numbers $N_C(\bm{v},\bm{a})$ 
of congruent classes of $4$-end catenoids 
with given $(\bm{v},\bm{a})$,
and a method to construct these surfaces.

\vspace{\baselineskip}

First we recall that 
if there is a $4$-end catenoid 
$x\colon{}\C\setminus\{q_1,q_2,q_3,q_4\}\to {\Real}^3$
such that $\nu(q_j)=v_j$ and $w(q_j)=a_j$ 
$(j=1,2,3,4)$, 
then, 
by \eqref{Reduction2}, 
there are nonzero complex numbers $b_1,b_2,b_3,b_4$ satisfying 
\[
    A
        \begin{pmatrix}
           b_1 \\
           b_2 \\
           b_3 \\
           b_4
        \end{pmatrix}
         = 
        \begin{pmatrix}
           0 \\
           0 \\
           0 \\
           0
        \end{pmatrix},
\]
with
\begin{equation}\label{mA}
     A := 
      \begin{pmatrix}
       0 &
       \displaystyle{\frac{\bar p_1p_2+\mathstrut{}1}{
            \mathstrut{}q_2-q_1}} &
       \displaystyle{\frac{\bar p_1p_3+\mathstrut{}1}{
            \mathstrut{}q_3-q_1}} &
       \displaystyle{\frac{\bar p_1p_4+\mathstrut{}1}{
            \mathstrut{}q_4-q_1}} \\
       \displaystyle{\frac{\bar p_2p_1+\mathstrut{}1}{
            \mathstrut{}q_1-q_2}} &
       0 &
       \displaystyle{\frac{\bar p_2p_3+\mathstrut{}1}{
            \mathstrut{}q_3-q_2}} &
       \displaystyle{\frac{\bar p_2p_4+\mathstrut{}1}{
            \mathstrut{}q_4-q_2}} \\
       \displaystyle{\frac{\bar p_3p_1+\mathstrut{}1}{
            \mathstrut{}q_1-q_3}} &
       \displaystyle{\frac{\bar p_3p_2+\mathstrut{}1}{
            \mathstrut{}q_2-q_3}} &
       0 &
       \displaystyle{\frac{\bar p_3p_4+\mathstrut{}1}{
            \mathstrut{}q_4-q_3}} \\
       \displaystyle{\frac{\bar p_4p_1+\mathstrut{}1}{
            \mathstrut{}q_1-q_4}} &
       \displaystyle{\frac{\bar p_4p_2+\mathstrut{}1}{
            \mathstrut{}q_2-q_4}} &
       \displaystyle{\frac{\bar p_4p_3+\mathstrut{}1}{
            \mathstrut{}q_3-q_4}} &
       0 
      \end{pmatrix},
\end{equation}
where $p_j:=\sigma(v_j)$ $(j=1,2,3,4)$ 
and $\sigma\colon{}{\bf S}^2\to\C$ 
is the stereographic projection 
from the north pole, 
that is, 
they satisfy the following identity: 
\[
     v_j = \left(
       \frac{2\Re(p_j)}{|p_j|^2+1},
       \frac{2\Im(p_j)}{|p_j|^2+1},
       \frac{|p_j|^2-1}{|p_j|^2+1}
       \right)
       \quad (j=1,\ldots,n). 
\]
Clearly, 
it holds that $\det A=0$.

To get upper estimates for $N_C(\bm{v},\bm{a})$, 
we discuss the rank of the matrix $A$.
We remark that 
$\rank A$ is invariant under 
both the conformal actions of the domain 
and the rigid motions in ${\Real}^3$.

When we consider $v$ of TYPE III, 
we may assume $p_4=q_4=\infty$ 
without loss of generality, 
and we can use the formula \eqref{Eq} 
in place of \eqref{Reduction2}.
Under this assumption, 
the matrix $A$ is given by 
\begin{equation}\tag{3.1$'$}\label{mA'}
     A = 
      \begin{pmatrix}
       0 &
       \displaystyle{\frac{\bar p_1p_2+\mathstrut{}1}{
            \mathstrut{}q_2-q_1}} &
       \displaystyle{\frac{\bar p_1p_3+\mathstrut{}1}{
            \mathstrut{}q_3-q_1}} &
       \bar p_1 \\
       \displaystyle{\frac{\bar p_2p_1+\mathstrut{}1}{
            \mathstrut{}q_1-q_2}} &
       0 &
       \displaystyle{\frac{\bar p_2p_3+\mathstrut{}1}{
            \mathstrut{}q_3-q_2}} &
       \bar p_2 \\
       \displaystyle{\frac{\bar p_3p_1+\mathstrut{}1}{
            \mathstrut{}q_1-q_3}} &
       \displaystyle{\frac{\bar p_3p_2+\mathstrut{}1}{
            \mathstrut{}q_2-q_3}} &
       0 &
       \bar p_3 \\
       -p_1 &
       -p_2 &
       -p_3 &
       0 
      \end{pmatrix}.
\end{equation}

\begin{Lemma}\label{lem:3-1}
    Let $\bm{v}=\{v_1,v_2,v_3,v_4\}$ be unit vectors of TYPE III\@.
    Then the number $N_q(\bm{v})$ of solutions $\bm{q}=\{q_1,q_2,q_3,q_4\}$ 
    of the equation $\det A=0$ is at most $4$ 
    up to M\"{o}bius transformations.
\end{Lemma}

\begin{pf}
    Set $p_j:=\sigma(v_j)$ as before.
    We may assume $p_4=q_4=\infty$ 
    without loss of generality, 
    and it follows from direct computation that 
    \[
        \det A
      = \frac{|p_1|^2|\bar p_2p_3+1|^2q_1^4+\cdots}
         {(q_1-q_2)^2(q_2-q_3)^2(q_3-q_1)^2}.
    \]
    Since we assume $\bm{v}$ is of TYPE III, 
    it is clear that 
    $|p_1|^2|\bar p_2p_3+1|^2\ne 0$.
    Therefore the number of $q_1$ satisfying $\det A=0$ 
    is at most $4$ 
    for any choice of $q_2$, $q_3$, 
    namely $N_q(\bm{v})\le 4$.
\end{pf}

\begin{Prop}\label{prop:3-2}
    $\rank A=3$ 
    if and only if $\bm{v}$ is of TYPE III\@.
\end{Prop}

\begin{pf}
    When $\bm{v}$ is of TYPE III, 
    we can easily see that 
    \[
         \det
         \begin{pmatrix}
              0 &
              \displaystyle{\frac{\bar p_2p_3+1}{q_3-q_2}} &
              \bar p_2 \\
              \displaystyle{\frac{\bar p_3p_2+1}{q_2-q_3}} &
              0 &
              \bar p_3 \\
              -p_2 &
              -p_3 &
              0 
         \end{pmatrix}
          = 
         \frac{\bar p_3p_2-\bar p_2p_3}
          {q_2-q_3}
        \ne 0.
    \]
    Hence we have $\rank A=3$
    for any $q_1,q_2,q_3$ satisfying $\det A=0$.
    
    On the other hand, 
    when $\bm{v}$ is of TYPE I or II, 
    by the remark above, 
    we may assume $p_1,p_2,p_3,p_4$ are real numbers.
    In this case, 
    since $A$ is skew-symmetric, 
    if $\det A=0$ and $A\ne 0$, 
    then we have $\rank A=2$.
    
    Now, 
    our assertion has been proved.
\end{pf}

\begin{Thm}\label{thm:3-3}
    For any $\bm{v}=\{v_1,v_2,v_3,v_4\}$ of TYPE III 
    and $\bm{a}=$ 
    $\{a_1$, $a_2$, $a_3$, $a_4\}$ satisfying 
    $\sum_{j=1}^4a_jv_j=0$, 
    $N_C(\bm{v},\bm{a})\le 4$. 
    Namely, 
    the number of $4$-end catenoids 
    with the same $(\bm{v}, \bm{a})$ of TYPE III 
    is at most $4$.
\end{Thm}

\begin{pf}
    From the proof of the proposition above, 
    $\dim\Ker A=1$ 
    for any $q_1,q_2,q_3$ satisfying $\det A=0$.
    Note that ${}^t(b_1,b_2,b_3,b_4)\in$ $\Ker A-\{0\}$.
    Moreover, 
    $(b_1,b_2,b_3,b_4)$ satisfies 
    $b_4(b_1+b_2+b_3)=a_4$ also.
    Hence, 
    if $(b_1,b_2,b_3,b_4)$ $=(\beta_1$, $\beta_2$, $\beta_3$, $\beta_4)$ 
    is a solution , 
    then the other solution is 
    $(b_1,b_2,b_3,b_4)=$ $-(\beta_1$, $\beta_2$, $\beta_3$, $\beta_4)$ 
    and both of these solutions give 
    the same Weierstrass data.
    Therefore, 
    it is clear that, 
    for any $q_1,q_2,q_3$ chosen above, 
    the number of $4$-end catenoids is at most $1$.
    Now we get the estimate 
    $N_C(\bm{v},\bm{a})\le N_q(\bm{v})\times 1\le 4$.
\end{pf}

\begin{Cor}\label{cor:3-4}
    Any $4$-end catenoid of TYPE III 
    is isolated in the sense of Rosenberg~\cite{rose}.
\end{Cor}

\begin{pf}
    Since $N_C(\bm{v},\bm{a})$ is finite 
    and any deformation moving flux 
    is not an $\epsilon$-$C^1$-variation, 
    our assertion is clear.
    \end{pf}

By solving the equation \eqref{Eq}
with $n=4$, 
$q_2=p_2$ and $q_3=p_3$ directly, 
we get the following method 
to construct $4$-end catenoids of TYPE III\@.

For given $(\bm{v},\bm{a})$, 
set $p_j:=\sigma(v_j)$ $(j=1,2,3,4)$ as before, 
and set 
\[
        A(t):=
        \begin{pmatrix}
           0 &
           \displaystyle{\frac{\bar p_1p_2+1}{p_2-t}} &
           \displaystyle{\frac{\bar p_1p_3+1}{p_3-t}} &
           \bar p_1 \\
           \displaystyle{\frac{\bar p_2p_1+1}{t-p_2}} &
           0 &
           \displaystyle{\frac{\bar p_2p_3+1}{p_3-p_2}} &
           \bar p_2 \\
           \displaystyle{\frac{\bar p_3p_1+1}{t-p_3}} &
           \displaystyle{\frac{\bar p_3p_2+1}{p_2-p_3}} &
           0 &
           \bar p_3 \\
           -p_1 &
           -p_2 &
           -p_3 &
           0 
        \end{pmatrix}
\]
(Remark that $A(q_1)=A|_{q_2=p_2,q_3=p_3}$),
\begin{align*}
    \Phi(t)&:=(p_2-p_3)^2(t-p_2)^2(t-p_3)^2 \det A(t)\\
           & = {\scriptstyle\det}\left(
               \begin{smallmatrix}
                   0 &
                   -(\bar p_1p_2+1)(t-p_3) &
                   -(\bar p_1p_3+1)(t-p_2) &
                   \bar p_1(t-p_2)(t-p_3) \\
                   (\bar p_2p_1+1)(p_2-p_3) &
                   0 &
                   -(\bar p_2p_3+1)(t-p_2) &
                   \bar p_2(p_2-p_3)(t-p_2) \\
                   (\bar p_3p_1+1)(p_2-p_3) &
                   (\bar p_3p_2+1)(t-p_3) &
                   0 &
                   \bar p_3(p_2-p_3)(t-p_3) \\
                   -p_1 &
                   -p_2 &
                   -p_3 &
                   0 
               \end{smallmatrix}\right)
           \\
           &= |p_1|^2|\bar p_2p_3+1|^2t^4+\cdots,
\end{align*}
and
\[
    B_1(t):=(\bar p_3p_2-\bar p_2p_3)(p_2-p_3)(t-p_2)(t-p_3).
\]
If $\det A(t)=0$ and $B_1(t)\ne 0$, 
then $\Ker A(t)$ is generated by 
${}^t(B_1(t)$, $B_2(t)$, $B_3(t)$, $B_4(t))$, 
where
\[
     \begin{pmatrix}
          B_2(t) \\
          B_3(t) \\
          B_4(t) 
     \end{pmatrix}
     := A'
     \begin{pmatrix}
          (\bar p_2p_1+1)(p_2-p_3)(t-p_3) \\
          (\bar p_3p_1+1)(p_2-p_3)(t-p_2) \\
          -p_1(t-p_2)(t-p_3) 
     \end{pmatrix},
\]
\[
     A' := 
     \begin{pmatrix}
          -|p_3|^2(p_2-p_3) &
          \bar p_2p_3(p_2-p_3) &
          \bar p_3(\bar p_2p_3+1)(p_2-p_3) \\
          \bar p_3p_2(p_2-p_3) &
          -|p_2|^2(p_2-p_3) &
          -\bar p_2(\bar p_3p_2+1)(p_2-p_3) \\
          p_3(\bar p_3p_2+1) &
          -p_2(\bar p_2p_3+1) &
          -|\bar p_2p_3+1|^2 
     \end{pmatrix}.
\]
(The matrix $A'$ is column-equivalent to 
the inverse of the $3\times 3$ submatrix 
which results by deleting 
the first row and column of the matrix $A(t)$.)

Note that $B_1(q_1)\ne 0$ holds 
for any solution $q_1$ 
of the equation $\Phi(t)=0$.
Indeed, 
if $B_1(q_1)=0$, 
then we have 
\[
     (\bar p_3p_2-\bar p_2p_3)(p_2-p_3)(q_1-p_2)(q_1-p_3)=0.
\]
However, 
since we assume $v$ is of TYPE III, 
\begin{alignat*}{2}
     \Phi(p_2) & = |p_3|^2|\bar p_1p_2+1|^2(p_2-p_3)^4 \,& \ne\, &0, \\
     \Phi(p_3) & = |p_2|^2|\bar p_1p_3+1|^2(p_2-p_3)^4 \,& \ne\, &0, 
\end{alignat*}
namely $q_1\ne p_2,p_3$.
Hence the equality above does not happen.
Assume 
$\prod_{k=2}^4B_k(q_1)\ne 0$ 
and 
$\sum_{j=1}^3B_j(q_1)\ne 0$.
Then, 
by straightforward calculation, 
we see that 
the solutions of the equation \eqref{Eq}
with $n=4$, 
$q_2=p_2$ and $q_3=p_3$ are given by 
\[
     \left\{
     \begin{aligned}
          q_1 &:\text{a solution of the equation}\quad \Phi(t)=0, \\
          q_2 &:= p_2, \quad q_3 := p_3, \quad  q_4 := \infty, \\
          b_1 &:= (\pm)B_1(q_1)
                 \sqrt{\displaystyle{
                 \frac{\mathstrut{}a_4}{
                    B_4(q_1)\sum_{j=1}^3B_j(q_1)}}}, \\
          b_j &:= 
            \displaystyle{\frac{B_j(q_1)}{B_1(q_1)}b_1} \qquad (j=2,3,4).
     \end{aligned}
     \right.
\]
If $q_1$ satisfies 
$\prod_{k=2}^4B_k(q_1)=0$ 
or 
$\sum_{j=1}^3B_j(q_1)=0$, 
then there are no solution $(\bm{q},\bm{b})$ 
of the equation \eqref{Eq} with such $q_1$, 
since we assume 
$a_j\ne 0$ $(j=1,2,3,4)$.
Set 
{\allowdisplaybreaks
\begin{align*}
\tilde P(z) & :=  p_1B_1(t)(z-p_2)(z-p_3)
                  +p_2B_2(t)(z-t)(z-p_3) \\
            &\hphantom{:=}\qquad    +p_3B_3(t)(z-t)(z-p_2) 
                    -B_4(t)(z-t)(z-p_2)(z-p_3), \\
\tilde Q(z) & :=  B_1(t)(z-p_2)(z-p_3)
                  +B_2(t)(z-t)(z-p_3)
                   +B_3(t)(z-t)(z-p_2), \\
\tilde R(z) & :=  (z-t)(z-p_2)(z-p_3). 
\end{align*}}
Note here that 
{\allowdisplaybreaks
    \begin{align*}
     \tilde P(z)|_{t=q_1} & = 
        \frac{B_1(q_1)}{b_1}P(z)|_{q_2=p_2,q_3=p_3}, \\
     \tilde Q(z)|_{t=q_1} & = 
        \frac{B_1(q_1)}{b_1}Q(z)|_{q_2=p_2,q_3=p_3}, \\
     \tilde R(z)|_{t=q_1} & = 
        R(z)|_{q_2=p_2,q_3=p_3}.
    \end{align*}}
Then it is easy to see that 
the $(\bm{q},\bm{b})$ above gives 
the following Weierstrass data 
of a $4$-end catenoid realizing $(\bm{v},\bm{a})$:
\begin{equation}\label{W4}
     g(z) = \left.
        \frac{\tilde P(z)}{\tilde Q(z)}
        \right|_{t=q_1},
     \qquad 
     \omega = -\frac{a_4}{B_4(q_1)\sum_{j=1}^3B_j(q_1)}
           \left(
           \left.
           \frac{\tilde Q(z)}{\tilde R(z)}
           \right|_{t=q_1}
           \right)^2dz.
\end{equation}
Let $\Psi(t)$ be the resultant of 
$\tilde P(z)$ and $\tilde Q(z)$.
The surface given by the data above 
has no branch points 
if and only if 
$q_1$ satisfies $\Psi(q_1)\ne 0$.

We can construct all of 
the $4$-end catenoids of TYPE III 
by this algorithm.

Now, 
we observe a typical 

\begin{Exa}\label{exa:3-5}
    Let $\zeta_3$ be a primitive root 
    of the equation $z^3=1$.
    For special values 
    $p_1=1/\sqrt{2}$, 
    $p_2=\zeta_3/\sqrt{2}$ and 
    $p_3=\zeta_3^2/\sqrt{2}$, 
    by direct computation, 
    we have 
    {\allowdisplaybreaks
    \begin{align*}
     \Phi(t) & =  \frac{3}{8}
                   \left(t-\frac{1}{\sqrt{2}}\right)^2(t+\sqrt{2})^2,\\
     B_1(t)  & =  \frac{3}{2\sqrt{2}}
                   \left(t^2+\frac{1}{\sqrt{2}}t+\frac{1}{2}\right),\\
     B_2(t)  & =  -\frac{3\zeta_3^2}{4\sqrt{2}}
                   \left(t^2-\frac{4\zeta_3+1}{\sqrt{2}}t-\zeta_3\right),\\
     B_3(t)  & =  -\frac{3\zeta_3}{4\sqrt{2}}
                   \left(t^2-\frac{4\zeta_3^2+1}{\sqrt{2}}t-\zeta_3^2\right),\\
     B_4(t)  & =  \frac{3}{4\sqrt{2}}
                   \left(t^2+\frac{1}{\sqrt{2}}t+2\right).
    \end{align*}}
    For one solution $1/\sqrt{2}$ 
    of the equation $\Phi(t)=0$, 
    we have
    {\allowdisplaybreaks
    \begin{align*}
         B_j\left(\frac{1}{\sqrt{2}}\right) 
              & =  \frac{9}{4\sqrt{2}} \qquad\qquad (j=1,2,3,4), \\
         \sum_{j=1}^3B_j\left(\frac{1}{\sqrt{2}}\right) 
              & =  \frac{27}{4\sqrt{2}} \ne 0, \\
         \tilde P(z)|_{t=1/\sqrt{2}} 
              & =  -\frac{9}{4\sqrt{2}}(z^3-\sqrt{2}), \\
         \tilde Q(z)|_{t=1/\sqrt{2}} 
              & =  \frac{27}{4\sqrt{2}}z^2, \\
         \tilde R(z)|_{t=1/\sqrt{2}} 
              & =  z^3-\frac{1}{2\sqrt{2}}, 
    \end{align*}}
    from which it follows that 
    $\Phi(1/\sqrt{2})\ne 0$.
    Now, 
    \eqref{W4} with these data gives 
    an Enneper-Weierstrass representation 
    of a tetrahedrally symmetric $4$-end cate\-noid.

    On the other hand, 
    for the other solution $-\sqrt{2}$, 
    we have
    {\allowdisplaybreaks
    \begin{align*}
         B_j(-\sqrt{2}) 
          & = \frac{9\zeta_3^{j-1}}{4\sqrt{2}} \qquad\qquad (j=1,2,3,4), \\
         \sum_{j=1}^3B_j(-\sqrt{2}) 
          & = 0, \\
         \tilde P(z)|_{t=-\sqrt{2}} 
          & =  -\frac{9}{4\sqrt{2}}
                \left(
                z^3+\frac{3}{\sqrt{2}}z^2+\frac{3}{2}z-\frac{1}{2\sqrt{2}}
                \right), \\
         \tilde Q(z)|_{t=-\sqrt{2}} 
             & =  -\frac{27}{8}\left(z+\frac{1}{\sqrt{2}}\right), \\
         \tilde R(z)|_{t=-\sqrt{2}} 
              & =  z^3+\frac{3}{\sqrt{2}}z^2+\frac{3}{2}z+\frac{1}{\sqrt{2}}, 
    \end{align*}}
    from which it follows that 
    $\Phi(-\sqrt{2})\ne 0$.
    Now, 
    the data 
    \[
         g(z) = \left.
            \frac{\tilde P(z)}{\tilde Q(z)}
            \right|_{t=q_1},
         \qquad 
         \omega = -\left(
            \left.
            \frac{\tilde Q(z)}{\tilde R(z)}
            \right|_{t=q_1}
            \right)^2dz.
    \]
    gives an Enneper-Weierstrass representation 
    of a complete minimal surface with $4$ flat ends.
    It is easy to see that 
    this surface is also 
    tetrahedrally symmetric 
    (cf. \cite{br}).

    From this consideration, 
    it is clear that 
    the tetrahedrally symmetric $4$-end catenoid 
    is unique up to homothety.
\end{Exa}

Now we will prove our second main theorem in Section~\ref{sec:1}.
Note here that 
$\Phi(t)$, 
$B_j(t)$ $(j=1,2,3,4)$ and $\Psi(t)$ 
are polynomials of $t$ 
whose coefficients are also polynomials of 
$p_1$, $p_2$, $p_3$, $\bar p_1$, $\bar p_2$ and $\bar p_3$.

\begin{Thm} \label{thm:3-6}
    For almost all pair $(\bm{v},\bm{a})$ 
    of unit vectors 
    $\bm{v}=\{v_1,v_2,v_3,v_4\}$ in ${\Real}^3$ 
    and nonzero real numbers 
    $\bm{a}=\{a_1,a_2,a_3,a_4\}$ 
    satisfying $\sum_{j=1}^4a_jv_j=0$, 
    there is a $4$-end catenoid 
    $x\colon{}\C\setminus\{q_1,q_2,q_3,q_4\}\to {\Real}^3$ 
    such that 
    $\nu(q_j)=v_j$ and the weight $w(q_j)=a_j$ $(j=1,2,3,4)$.
\end{Thm}

\begin{pf}
    Set $p_j:=\sigma(v_j)$ and 
    assume $p_4=q_4=\infty$ as before.
    Then $\bm{v}$ is of TYPE III 
    if and only if 
    \begin{align*}
         D & :=  (\bar p_2p_3-p_2\bar p_3)
            (\bar p_3p_1-p_3\bar p_1)
            (\bar p_1p_2-p_1\bar p_2)
            \times\{(|p_1|^2-1)(\bar p_2p_3-p_2\bar p_3)\\
           & \hphantom{:=}\qquad
            +(|p_2|^2-1)(\bar p_3p_1-p_3\bar p_1)
            +(|p_3|^2-1)(\bar p_1p_2-p_1\bar p_2)\}\\
           & \ne  0.
    \end{align*}
    ($D$ can be obtained from the multiplication 
    of the $3\times 3$ minor determinants 
    of the $3\times 4$ matrix $(v_1,v_2,v_3,v_4)$, 
    where $v_4={}^t(0,0,1)$).
    In this case, 
    what we have only to do is 
    to show the existence of a solution $q_1$ 
    of the equation $\Phi(t)=0$ satisfying 
    $\prod_{k=2}^4B_k(q_1)\ne 0$, 
    $\sum_{j=1}^3B_j(q_1)\ne 0$ and $\Psi(q_1)\ne 0$.
    Let $\Omega(t)$ be the remainder upon division of 
    $(\Psi(t)\prod_{k=2}^4B_k(t)\sum_{j=1}^3B_j(t))^4$ 
    by $\Phi(t)$.
    Clearly, 
    $\Omega(t)\not\equiv 0$ 
    if and only if 
    there exists at least one solution $q_1$ 
    of the equation $\Phi(t)=0$ satisfying 
    $\prod_{k=2}^4B_k(q_1)\ne 0$, 
    $\sum_{j=1}^3B_j(q_1)\ne 0$ and $\Psi(q_1)\ne 0$.

    Now, 
    since there exists 
    the  tetrahedrally symmetric $4$-end catenoid, 
    it is clear that 
    at least one coefficient $\Omega_0$ of $\Omega(t)$ 
    does not vanish 
    for special values 
    $p_1=1/\sqrt{2}$, 
    $p_2=\zeta_3/\sqrt{2}$ and 
    $p_3=\zeta_3^2/\sqrt{2}$ 
    (see Example \ref{exa:3-5}).
    Remark here that 
    each coefficient of $\Omega(t)$ is 
    a rational function of 
    $p_1,p_2,p_3,\bar p_1,\bar p_2,\bar p_3$, 
    i.e.~it is real analytic.
    Hence, 
    we see that $\Omega(t)\not\equiv 0$ 
    for almost all $p_1,p_2,p_3$. 
    Now we have proved that, 
    for any $p_1, p_2, p_3, p_4(=\infty)$ 
    satisfying 
    $(\Omega_0D)(p_1,p_2,p_3,\bar p_1,\bar p_2,\bar p_3)\ne 0$ 
    and $\bm{a}=\{a_1,a_2,a_3,a_4\}$ such that 
    $\sum_{j=1}^4a_j\sigma^{-1}(p_j)=0$, 
    there is at least one $4$-end catenoid 
    $x\colon{}\C\setminus\{q_1,q_2,q_3,q_4\}\to {\Real}^3$ 
    such that $\nu(q_j)=\sigma^{-1}(p_j)$ and $w(q_j)=a_j$.
\end{pf}

We know that 
we cannot remove \lq\lq almost".
Indeed, 
there are additional obstructions 
for the existence of $4$-(or $n$-)end catenoids.
We will consider these 
in the following section.

Before concluding this section, 
let us observe a significant example 
which includes 
a deformation from 
the Jorge-Meeks $4$-noid 
to the tetrahedrally symmetric $4$-end catenoid.
Moreover, we will remark there that
the symmetry of the flux data of an $4$-end catenoid 
does not always imply the symmetry of the surface.

\begin{Exa}\label{exa:3-7}
We will now try to construct a $4$-end catenoid 
with the following data;

\begin{equation}\label{datap}
    \left\{
    \begin{aligned}
     & p_1=p, \quad p_2=-p, \quad p_3=p^{-1}i, \quad p_4=-p^{-1}i, \\
     & a_1=a_2=a_3=a_4=1, 
    \end{aligned}
    \right.
\end{equation}
where $p$ is a positive real number. 
The corresponding flux data $(\bm{v}, \bm{a})$ 
is invariant under the action of the group 
\[
    \left\langle 
    \left(
    \begin{array}{rrr}
     -1 &  0 &  0 \\
      0 & -1 &  0 \\
      0 &  0 &  1 
    \end{array}
    \right),
    \left(
    \begin{array}{rrr}
      0 &  1 &  0 \\
      1 &  0 &  0 \\
      0 &  0 & -1 
    \end{array}
    \right),
    \left(
    \begin{array}{rrr}
     -1 &  0 &  0 \\
      0 &  1 &  0 \\
      0 &  0 &  1 
    \end{array}
    \right)
    \right\rangle 
    \subset O(3). 
\]
We may put
\[
    q_1=q, \quad q_2=-q, \quad q_3=q^{-1}i, \quad q_4=-q^{-1}i, \\
\]
for a nonzero complex number $q$ $(q^4\ne -1)$ 
without loss of generality, 
and the determinant of the matrix $A$ in (3.1) 
is computed as 
\[
    \det A = \left[\left\{\frac{2q(q^2-1)}{q^4+1}\right\}^2 
                   -\left(\frac{p^2-1}{2p}\right)^2 
             \right]
             \cdot 
             \left[\left\{\frac{2q(q^2+1)}{q^4+1}\right\}^2 
                   +\left(\frac{p^2-1}{2p}\right)^2 
             \right].
\]
Thus $\det A=0$ if and only if 
one of the following four equations is satisfied: 
\begin{gather}
    \frac{2q(q^2-1)}{q^4+1} = \pm\frac{p^2-1}{2p}; 
    \tag{$3.4_\pm$} \\ 
    \frac{2q(q^2+1)}{q^4+1} = \pm\frac{p^2-1}{2p}i. 
    \tag{$3.5_\pm$} 
\end{gather}
Note here that 
$q=\alpha$ is a solution of $(3.4_+)$ 
if and only if 
$q=-\alpha$ (resp. $q=\pm\alpha^{-1}i$) 
is a solution of $(3.4_-)$ (resp. $(3.5_\pm)$), 
and that the corresponding solutions $(\bm{q}, \bm{b})$ 
of \eqref{Reduction2} 
give the Weierstrass data of the surfaces 
congruent to each other. 
Therefore we may only consider $(3.4_+)$. 
It is easy to see that 
the equation $(3.4_+)$ has a real solution $q$ 
if and only if $c^{-1}\le p\le c$, 
where  $c:=(\sqrt{6}+\sqrt{2})/2$. 

Let $q$ be a solution of the equation $(3.4_+)$. 
Now, 
since our data \eqref{datap} is of TYPE III, 
rank$A$ must be $3$ 
(except for the case $p=1$). 
It is clear that the nonzero vector 
${}^t(q,q,p,p)\in$Ker$A$, 
and hence it spans Ker$A$. 
Therefore we can set
\[
     b_1=b_2=rq,\qquad b_3=b_4=rp, 
\]
and, 
by \eqref{Reduction2}, 
we get 
\[
     r^2 = \displaystyle{
           \frac{q^{4}+1}
                {q\{p(q^{4}+1)+2q(p^2q^{2}+1)\}}
                        }. 
\]

Define a surface 
$x_{q}\colon{}\C\setminus\{q,-q,q^{-1}i,-q^{-1}i\} \to {\Real}^3$ 
by these data, 
whose Weierstrass data is given by 
\begin{gather}
    g(z) := \frac{(pq^2-q^{-1})z^2+(p+q)} 
            {z\{(p+q)z^2-(pq^2-q^{-1})\}}, 
    \tag{3.6} 
    \\
    \omega := -r^2\left[\frac{2z\{(p+q)z^2-(pq^2-q^{-1})\}} 
                             {(z^2-q^2)(z^2+q^{-2})} 
                  \right]^2dz. 
    \tag{3.7} 
\end{gather}
Then $x_q$ is a branched conformal minimal immersion 
such that 
\[
    \left\{\begin{array}{l}
        g(q)=p, \quad 
        g(-q)=-p, \quad 
        g(q^{-1}i)=p^{-1}i, \quad 
        g(-q^{-1}i)=-p^{-1}i, \\ 
        w(q)=w(-q)=w(q^{-1}i)=w(-q^{-1}i)=1, 
           \end{array}
    \right.
\]
except for the case when 
$q=-c^{\pm 1}$ (i.e. $r=\infty$). 
It is easy to see that 
\begin{enumerate}
    \item[(i)] 
      $x_q$ is branched if and only if $q=-1$ $(p=1)$;
    \item[(ii)] 
      The normalized surface $\tilde x_q$ 
      defined by the same $g$ as $x_q$ and $\tilde\omega:=\omega/r^2$ 
      has $4$ flat ends (i.e. $a_j=0$ for any $j$) 
      if and only if $q=-c^{\pm 1}$ and $p=c^{\mp 1}$. 
    \item[(iii)] 
      In the other cases, 
      $x_q$ is a non-branched $4$-end catenoid. 
      In particular, 
      when $c^{-1}\le p\le c$, 
      any solution $q$ of $(3.4_+)$ is real, 
      and $x_q$ has the same symmetry as its flux data (\ref{datap}). 
      On the other hand, 
      when $0<p<c^{-1}$ or $c<p$, 
      any solution $q$ of $(3.4_+)$ 
      can take neither a real nor a purely imaginary value, 
      and the isometry group of $x_q$ is $\Z_2\times\Z_2$ 
      which is smaller than that of its flux data (\ref{datap}). 
\end{enumerate}
In particular, 
\begin{enumerate}\setcounter{enumi}{3}
    \item[(iv)] 
      $x_1$ is the Jorge-Meeks $4$-noid; 
    \item[(v)] 
      $x_c$ is the tetrahedrally symmetric $4$-end catenoid. 
\end{enumerate}

By (iii), (iv) and (v), 
we see that 
the family $\{x_{q};\, 1\le q\le c\}$ 
gives a deformation from the Jorge-Meeks $4$-noid 
to the tetrahedrally  symmetric one. 
In other words, 
the Jorge-Meeks $4$-noid  
and the tetrahedrally  symmetric one 
are included in the same connected component 
of the moduli of $4$-end  catenoids. 

Furthermore, 
we can check that, 
for any real number $p$ such that $1<p<c$, 
there are four real numbers 
$q^{(1)},q^{(2)},q^{(3)}$ and $q^{(4)}$ 
satisfying 
\[
        q^{(1)}>q^{(2)}>1>0
         >q^{(3)}(=-\frac{1}{q^{(1)}})
         >q^{(4)}(=-\frac{1}{q^{(2)}})>-1
\] 
and $(3.4_+)$ with 
$q=q^{(\ell)}$ ($\ell=1,2,3,4$).
Hence there are four $4$-end catenoids 
$x_{q^{(1)}}$, $x_{q^{(2)}}$, $x_{q^{(3)}}$ and $x_{q^{(4)}}$ 
which have the same flux data.
It can be easily observed that 
these four surfaces are not congruent to each other.
This concludes that 
our estimate $N_C(\bm{v},\bm{a})\le 4$ 
in Theorem~\ref{thm:3-3} is sharp.

These situations are demonstrated in Figure~\ref{fig:graph}.
Figure \ref{fig:4noid2} 
shows the image of $x_q$ for various values of $q$ 
with the same flux.
\end{Exa}

\begin{figure}[thb]
        \vspace{7cm}
        \caption{Example \ref{exa:3-7}.}
        \label{fig:graph}
\end{figure}

\begin{figure}[thb]
        \hspace*{\fill}
        \begin{minipage}{6cm}
        \footnotesize
            \vspace{5cm}
            \begin{center}
                (a) $p=1.2$, $q\approx 1.0976$
            \end{center}
        \end{minipage}
            \hfill
        \begin{minipage}{6cm}
        \footnotesize
        \vspace{5cm}
            \begin{center}
                (b) $p=1.2$, $q\approx 10.815$
            \end{center}
        \end{minipage}
        \hfill\\
        \hspace*{\fill}
        \begin{minipage}{6cm}
        \footnotesize
        \vspace{5cm}
            \begin{center}
                (c) $p=1.2$, $q\approx -0.91078$
            \end{center}
        \end{minipage}
            \hfill
        \begin{minipage}{6cm}
        \footnotesize
        \vspace{5cm}
            \begin{center}
                (d) $p=1.2$, $q\approx -0.09246$
            \end{center}
        \end{minipage}
        \hfill
        \caption{Non-congruent four surfaces with the same flux.}
        \label{fig:4noid2}
\end{figure}


\section{$4$-end catenoids of special type and additional obstructions} 
\label{sec:4} 

In this section, 
we consider the cases of TYPEs I and II\@.
First, 
we will prove that 
the same assertions as in Theorem \ref{thm:3-3} and Corollary \ref{cor:3-4}
hold  also in the case of TYPE II\@.

Let $N_C(\bm{v},\bm{a})$, 
$N_q(\bm{v})$ and $A$ 
be as in the previous section (see Lemma \ref{lem:3-1} etc.).

\begin{Lemma}\label{lem:4-1}
    Let $\bm{v}=\{v_1,v_2,v_3,v_4\}$ be 
    a family of unit vectors of TYPE II\@.
    Then $N_q(\bm{v})\le 2$.
\end{Lemma}

\begin{pf}
    Set $p_j\colon{}=\sigma(v_j)$ $(j=1,2,3,4)$ as before, 
    where $\sigma$ is the stereographic projection.
    We may assume $p_1,p_2,p_3,p_4$ are real numbers 
    without loss of generality.
    Set $p_{jk}:=p_jp_k+1$ for convenience.
    It follows from direct computation that
    \[
         \det A 
          = \left\{
            \displaystyle{  
                 \frac{p_{12}p_{34}}
                    {(q_1-q_2)(q_3-q_4)}
                -\frac{p_{13}p_{24}}
                      {(q_1-q_3)(q_2-q_4)}
                +\frac{p_{14}p_{23}}
                      {(q_1-q_4)(q_2-q_3)}
            }
            \right\}^2.
    \]
    We may also assume 
    $q_1=-q_2\ne 0,\pm 1$ 
    and 
    $q_3=-q_4=1$ 
    without loss of generality, 
    and we have
    \[
         \det A = 
             \frac{q_1\Phi_{II}(q_1+q_1^{-1})}
                  {4(q_1-1)^2(q_1+1)^2},
    \]
    where
    \[
         \Phi_{II}(t) = p_{12}p_{34}t^2
               +4(p_{13}p_{24}-p_{14}p_{23})t
               +8(p_{13}p_{24}+p_{14}p_{23})-4p_{12}p_{34}.
    \]
    Since we assume $\bm{v}$ is of TYPE II, 
    it is clear that 
    at least one of the coefficients of $\Phi_{II}(t)$ 
    does not vanish.
    Therefore 
    the number of $q_1+q_1^{-1}$ satisfying det$A=0$ 
    is at most $2$.
    Since $\C\setminus\{\pm q_1,\pm 1\}$ 
    and $\C\setminus\{\pm q_1^{-1},\pm 1\}$ 
    are conformal to each other 
    by a M\"{o}bius transformation, 
    we get $N_q(\bm{v})\le 2$.
\end{pf}

\begin{Thm}\label{thm:4-2}
    For any $\bm{v}=\{v_1,v_2,v_3,v_4\}$ of TYPE II 
    and $\bm{a}=\{a_1,a_2,a_3,a_4\}$ satisfying 
    $\sum_{j=1}^4a_jv_j=0$, 
    $N_C(\bm{v},\bm{a})\le 4$. 
    Namely, 
    the number of $4$-end catenoids 
    with the same $(\bm{v},\bm{a})$ of TYPE II 
    is at most $4$.
\end{Thm}

\begin{pf}
    From the proof of Proposition \ref{prop:3-2}, 
    $\dim\Ker A=2$ 
    for any $q_1,q_2,q_3,q_4$ satisfying $\det A=0$.
    Note that ${}^t(b_1,b_2,b_3,b_4)\in$ $\Ker A-\{0\}$.
    We may assume $p_{12}p_{34}(p_1-p_2)(p_3-p_4)\ne 0$ 
    without loss of generality.
    By putting $q_1=-q_2$ and $q_3=-q_4=1$ into \eqref{Reduction2}, 
    we have
    \begin{equation}\label{bb}
        \left\{
        \begin{aligned}
          b_3 &= \displaystyle{\frac 2{p_{34}}}
            \left(
            \displaystyle{
            -\frac{p_{14}}{q_1+1}b_1
            +\frac{p_{24}}{q_1-1}b_2
            }
            \right) \\
          b_4 &= -\displaystyle{\frac 2{p_{34}}}
                \left(
            \displaystyle{
                -\frac{p_{13}}{q_1-1}b_1
            +\frac{p_{23}}{q_1+1}b_2
                }
            \right), 
        \end{aligned}\right.
    \end{equation}
    \begin{equation}\label{bbaa}
        \left\{
        \begin{aligned}
          b_3\left(
          \displaystyle{
          \frac{p_1-p_3}{q_1-1}b_1
            -\frac{p_2-p_3}{q_1+1}b_2
            -\frac{p_4-p_3}2b_4}
            \right)
            &= a_3 \\
          b_4\left(
            \displaystyle{
            \frac{p_1-p_4}{q_1+1}b_1
            -\frac{p_2-p_4}{q_1-1}b_2
            +\frac{p_3-p_4}2b_3}
            \right)
            &= a_4.
         \end{aligned}\right.
    \end{equation}
    Putting \eqref{bb} into \eqref{bbaa}, 
    we get
\footnotesize
    \begin{equation} \label{pp}
         \left\{
         \begin{array}{r}
            \dfrac{p_3^2+1}{a_3}
          \left[
           -\dfrac{p_{14}(p_1-p_4)}
          {(q_1+1)(q_1-1)}
            b_1^2
           +\left\{
            \dfrac{p_{24}(p_1-p_4)}
                  {(q_1-1)^2}
           +\dfrac{p_{14}(p_2-p_4)}
                  {(q_1+1)^2}
            \right\}
            b_1b_2 
          \right. \\
          \left.
           -\dfrac{p_{24}(p_2-p_4)}
                  {(q_1+1)(q_1-1)}
            b_2^2
          \right]
          = \dfrac{p_{34}^2}2 \\
           -\dfrac{p_4^2+1}{a_4}
          \left[
           -\dfrac{p_{13}(p_1-p_3)}
                  {(q_1+1)(q_1-1)}
            b_1^2
           +\left\{
            \dfrac{p_{13}(p_2-p_3)}
                  {(q_1-1)^2}
           +\dfrac{p_{23}(p_1-p_3)}
                  {(q_1+1)^2}
            \right\}
           b_1  b_2 
          \right. \\
           \qquad\qquad
          \left.
        -\dfrac{p_{23}(p_2-p_3)}
          {(q_1+1)(q_1-1)}
            b_2^2
            \right]
            = \dfrac{p_{34}^2}2.
        \end{array}
        \right.
    \end{equation}
\normalsize
    If this equation has more than $4$ solutions $(b_1,b_2)$ 
    such that $b_1b_2\ne 0$, 
    then it holds that 
    \[
        \left\{
        \begin{aligned}
          \dfrac{p_3^2+1}{a_3}
            \left\{
            -\dfrac{p_{14}(p_1-p_4)}
                {(q_1+1)(q_1-1)}
            \right\}
            &= 
            -\dfrac{p_4^2+1}{a_4}
            \left\{
            -\dfrac{p_{13}(p_1-p_3)}
            {(q_1+1)(q_1-1)}
        \right\}
        \\
        \frac{p_3^2+1}{a_3}
        \left\{
            -\dfrac{p_{24}(p_2-p_4)}
            {(q_1+1)(q_1-1)}
            \right\}
            &= 
         -\dfrac{p_4^2+1}{a_4}
        \left\{
            -\dfrac{p_{23}(p_2-p_3)}
                {(q_1+1)(q_1-1)}
        \right\}, 
        \end{aligned}\right.
    \]
    namely
    \[
        \begin{pmatrix}
          p_{13}(p_1-p_3) & 
          p_{14}(p_1-p_4) \\ 
          p_{23}(p_2-p_3) & 
          p_{24}(p_2-p_4) 
        \end{pmatrix}
        \begin{pmatrix}
          \displaystyle{\frac{a_3}{p_3^2+1}} \\
          \displaystyle{\frac{a_4}{p_4^2+1}} 
        \end{pmatrix}
         = 
        \begin{pmatrix}
          0 \\ 0 
        \end{pmatrix}.
    \]
    Hence 
    \[
         0 = \det 
         \begin{pmatrix}
          p_{13}(p_1-p_3) & 
          p_{14}(p_1-p_4) \\ 
          p_{23}(p_2-p_3) & 
          p_{24}(p_2-p_4) 
         \end{pmatrix}
             = p_{12}p_{34}(p_1-p_2)(p_3-p_4).
    \]
    This contradicts our assumption.
    Therefore, 
    the equation \eqref{pp} has at most $4$ solutions $(b_1,b_2)$
    such that $b_1b_2\ne 0$.
    Note that if $(b_1,b_2)=(\beta_1,\beta_2)$ is a solution, 
    then $(b_1,b_2)=-(\beta_1,\beta_2)$ is also a solution 
    and both of these solutions give the same Weierstrass data.
    Therefore it is clear that, 
    for any conformal class chosen above, 
    the number of $4$-end catenoids realizing the given $(\bm{v},\bm{a})$ 
    is at most $2$.
    Now we get the estimate $N_C(\bm{v},\bm{a})\le N_q(\bm{v})\times 2 \le 4$.
\end{pf}

\begin{Cor}\label{cor:4-3}
    Any $4$-end catenoid of TYPE II 
    is isolated in the sense of Rosenberg~\cite{rose}.
\end{Cor}

\begin{pf}
    By the same reason as in the proof of Corollary \ref{cor:3-4}, 
    our assertion is clear.
\end{pf}

\begin{Rmk}
    In Lemma~\ref{lem:4-1},
    $N_q(\bm{v})\le 1$ 
    if one of the following condition holds:
    \begin{enumerate}
        \item $p_{jk}=0$, i.e. 
              $-v_j=v_k$ for some $j\ne k$;
        \item $D:={p_{12}}^2{p_{34}}^2+{p_{13}}^2{p_{24}}^2
                  +{p_{14}}^2{p_{23}}^2
                  -2(p_{13}p_{24}p_{14}p_{23}+p_{12}p_{34}p_{14}p_{23}
                  +p_{12}p_{34}p_{13}p_{24})=0.$
    \end{enumerate}
    Indeed,
    it holds that 
\footnotesize
    \[
        \Phi_{II}(q_1+q_1^{-1}) = 
        \begin{cases}
            4(p_{13}p_{24}-p_{14}p_{23})(q_1+q_1^{-1})
            +8(p_{13}p_{24}+p_{14}p_{23}) 
            \quad \mbox{if} \quad p_{12}p_{34}=0 \\
            \left\{p_{12}p_{34}(q_1+q_1^{-1})
            -4p_{14}p_{23}+2p_{12}p_{34}\right\}
            \displaystyle{\frac{(q_1-1)^2}{q_1}} 
            \quad \mbox{if} \quad p_{13}p_{24}=0 \\
            \left\{p_{12}p_{34}(q_1+q_1^{-1})
            +4p_{13}p_{24}-2p_{12}p_{34}\right\}
            \displaystyle{\frac{(q_1+1)^2}{q_1}} 
            \quad \mbox{if} \quad p_{14}p_{23}=0 \\
            \left\{p_{12}p_{34}(q_1+q_1^{-1})
            +2p_{13}p_{24}-2p_{14}p_{23}\right\}^2
            \displaystyle{\frac{1}{p_{12}p_{34}}} 
            \quad \mbox{if} \quad D=0.
        \end{cases}
    \]
\normalsize
    Therefore,
    by the proof of Theorem~\ref{thm:4-2},
    if either (1) or (2) holds,
    then we get the estimate 
    $N_C(\bm{v},\bm{a})\le N_q(\bm{v})\times 2\le 2$.
\end{Rmk}

Moreover, 
we can find an additional obstruction 
for the existence of $4$-end catenoids of TYPE II.

\begin{Thm}\label{thm:4-5}
    There are no $4$-end catenoids 
    $x\colon{}\C\setminus\{q_1,q_2,q_3,q_4\}\to {\Real}^3$
    such that $\nu(q_j)=v_j$ and $w(q_j)=a_j$ $(j=1,2,3,4)$ 
    if $-v_1=v_2$ and $v_3=v_4\ne\pm v_1$.
\end{Thm}

\begin{pf} 
    Set $p_j:=\sigma(v_j)$ as before.
    When our assumption holds, 
    we may assume $p_1,p_2,p_3,p_4$ 
    are nonzero real numbers satisfying 
    $p_1p_2+1=0$ and $p_3=p_4\ne p_1,-p_1^{-1}$ 
    without loss of generality.
    Suppose there exists a $4$-end catenoid 
    $x\colon{}\C\setminus\{q_1,q_2,q_3,q_4\}\to {\Real}^3$
    such that $\nu(q_j)=v_j$ and $w(q_j)=a_j$ $(j=1,2,3,4)$ 
    for some nonzero real numbers $a_1,\ldots,a_n$ 
    satisfying $a_1=a_2$ and $a_3+a_4=0$.
    Then it follows from direct computation that 
    \[
         0 = \det A 
       = \left\{
         \frac{(p_1p_3+1)(p_1-p_3)(q_1-q_2)(q_3-q_4)}
          {p_1(q_1-q_3)(q_2-q_3)(q_1-q_4)(q_2-q_4)}
        \right\}^2.
    \]
    This contradicts our assumption $p_3\ne p_1,-p_1^{-1}$
\end{pf}

Next, 
we consider the case of TYPE I, 
namely, 
the case when all of the ends are parallel.
In this case, 
there exist additional obstructions for an arbitrary $n$. 
For instance, 
since the degree of the Gauss map must be less than $n$ 
for any $n$-end catenoid (see Section \ref{sec:2}), 
the flux data 
\begin{equation}\label{obs3}
    v_1=\cdots=v_n 
\end{equation}
cannot be realized by any $n$-end catenoid. 
Moreover, 
the following obstructions can be also 
found by using (\ref{Reduction2}). 
\begin{eqnarray}
 & -v_1=-v_2=v_3=\cdots=v_n \label{obs1} \\
 & -v_1=v_2=\cdots=v_n, \quad 
    \sum_{j=2}^{n-1}\sum_{k=j+1}^na_ja_k\ne 0. \label{obs2} 
\end{eqnarray}
The obstruction \eqref{obs2} follows 
from the compatibility condition in Perez~\cite{per}. 
Conversely, 
in the exceptional case of the obstruction \eqref{obs2}, 
we have the following 

\begin{Lemma}\label{lem:4-7}
    Let $\bm{v}=\{v_1,\ldots,v_n\}$ be 
    a family of unit vectors in ${\Real}^3$, 
    and $\bm{a}=\{a_1,\ldots,a_n\}$ 
    a set of nonzero real numbers satisfying 
    $-v_1=v_2=\cdots=v_n$, 
    $a_1=\sum_{j=2}^na_j$ 
    and $\sum_{j=2}^{n-1}\sum_{k=j+1}^na_ja_k=0$.
    Then there exists an n-end catenoid 
    $x\colon{}\C\setminus\{q_1,\ldots,q_n\}\to {\Real}^3$ 
    such that $\nu(q_j)=v_j$ and $w(q_j)=a_j$ $(j=1,\ldots,n)$ 
    if and only if there are complex numbers 
    $q_1,\ldots,q_n$ satisfying 
    \begin{equation}\label{Eqf}
         \Condsum{k=2}{k\ne j}^n\frac{a_k}{q_k-q_j} = 0 
         \qquad (j=2,\ldots,n).
    \end{equation}
\end{Lemma}

\begin{pf}
    Set $p_j:=\sigma(v_j)$ as before. 
    We may assume $p_1=q_1=\infty$ 
    and $p_2=\cdots=p_n=0$ 
    without loss of generality.
    Suppose there exists an $n$-end catenoid 
    $x\colon{}\C\setminus\{q_1,\ldots,q_n\}\to {\Real}^3$
    such that $\nu(q_j)=v_j$ and $w(q_j)=a_j$ $(j=1,\ldots,n)$.
    Then, 
    by \eqref{Eq} 
    with the assumption $p_1=q_1=\infty$ 
    in place of $p_n=q_n=\infty$, 
    it holds that 
    \[
         b_jb_1=a_j, 
         \quad 
         \Condsum{k=2}{k\ne j}^n\frac{b_k}{q_k-q_j}=0
         \qquad
         (j=2,\ldots,n).
    \]
    Hence we have
    \[
         0 = b_1\Condsum{k=2}{k\ne j}^n\frac{b_k}{q_k-q_j}
           = \Condsum{k=2}{k\ne j}^n\frac{a_k}{q_k-q_j}
         \qquad
         (j=2,\ldots,n).
    \]
    Conversely, 
    suppose there are $n$ complex numbers $q_1,\ldots,q_n$ 
    satisfying \eqref{Eqf}.
    Put
    \[
         f(z):=\sum_{j=2}^n\frac{a_j}{z-q_j},
    \]
    and set, 
    for any nonzero complex number $t$,
    \begin{equation}\label{Wf}
         g_t(z) := -\frac 1{tf(z)},
         \quad
         \omega_t := -t(f(z))^2dz.
    \end{equation}
    Then, 
    for any $t$, 
    the surface 
    $x_t\colon{}\C\setminus\{q_1,\ldots,q_n\}\to {\Real}^3$
    represented by these data 
    is an $n$-end catenoid 
    such that $\nu(q_j)=v_j$, 
    $w(q_j)=a_j$  $(j=1,\ldots,n)$ 
    and the induced metric 
    \[
         ds_t^2 = \frac{(|tf|^2+1)^2}{|t|^2}|dz|^2
    \]
    is non-degenerate.
\end{pf}

By the proof of Lemma \ref{lem:4-7},
if there are $n$ complex numbers $q_1,\ldots,q_n$ 
satisfying \eqref{Eqf}, 
then there exists a 1-parameter family 
$\{x_t;\,t\in{\C}\setminus\{0\}\}$ 
of $n$-end catenoids with the same $(\bm{v}, \bm{a})$.
However, 
when $|t|=|t'|$, 
$x_t$ can be transformed to $x_{t'}$ 
by a certain rotation around the $x_3$-axis.
On the other hand, 
in the case $n>2$, 
when $|t|\ne|t'|$,
clearly $x_t$ and $x_{t'}$ are 
not isometric to each other.
Hence, 
the family $\{x_t;\,t>0\}$ is 
a non-trivial deformation.
 \par

Note that this deformation is 
an example of the deformation described 
in Lopez-Ros~\cite{lr}.
It is clear that $x_t$ is deformable 
in the sense of Rosenberg~\cite{rose}.

By solving the equation \eqref{Eqf} with $n=4$ and $5$, 
we have the following examples 
which completes the classification 
of at most $5$-end catenoids of parallel ends.
Indeed, 
by virtue of the conditions \eqref{obs3}, \eqref{obs1} and \eqref{obs2}, 
it is clear that 
any $n$-end catenoid of TYPE I with $n\le 5$ 
coincides with the ($2$-end) catenoid 
or one of the surfaces in Examples \ref{exa:4-8} and \ref{exa:4-9}. 

\begin{Exa}\label{exa:4-8}
    Let $\bm{v}=\{v_1,v_2,v_3,v_4\}$ 
    be a family of unit vectors in ${\Real}^3$ 
    satisfying $-v_1=v_2=v_3=v_4$.
    For any set 
    $\bm{a}=\{a_1,a_2,a_3,a_4\}$ 
    of nonzero real numbers 
    satisfying $a_1=a_2+a_3+a_4$ and 
    $a_2a_3+a_2a_4+a_3a_4=0$, 
    there exist a unique 1-parameter family 
    $\{x_t\colon{}\C\setminus\{q_1,q_2,q_3,q_4\}\to {\Real}^3
        ;\,t>0\}$ 
    of $4$-end catenoids  
    such that $\nu(q_j)=v_j$ and $w(q_j)=a_j$ $(j=1,2,3,4)$.
    Indeed their representations are given by \eqref{Wf} with 
    \[
        f(z) := \frac{a_2}z
         + \frac{a_3}{z-1}
         + \frac{a_4}
                {z+\displaystyle{\frac{a_4}{a_3}}},
    \]
    up to rigid motion in ${\Real}^3$.
    Figure \ref{fig:5-9} shows 
    the image of $x_{t}$ for various value of $t$, 
    when $a_2=-1$ and $a_3=a_4=2$.
\end{Exa}

    \begin{figure}[thb]
        \begin{minipage}{45mm}
        \vspace{5cm}
            \footnotesize
            \begin{center}
                (a) $t=0.5$
            \end{center}
        \end{minipage}
        \hfill
        \begin{minipage}{45mm}
        \vspace{5cm}
            \footnotesize
            \begin{center}
                (b) $t=1$
            \end{center}
        \end{minipage}
        \hfill
        \begin{minipage}{45mm}
        \vspace{5cm}
            \footnotesize
            \begin{center}
                (c) $t=2$
            \end{center}
        \end{minipage}
    \caption{Example~\ref{exa:4-8} for $a_2=-1$ and $a_3=a_4=2$}
    \label{fig:5-9}
    \end{figure}

\begin{Exa}\label{exa:4-9}
    Let $\bm{v}=\{v_1,v_2,v_3,v_4,v_5\}$ 
    be a family of unit vectors in ${\Real}^3$ 
    satisfying $-v_1=v_2=v_3=v_4=v_5$.
    For any set 
    $\bm{a}=\{a_1,a_2,a_3,a_4,a_5\}$ 
    of nonzero real numbers 
    satisfying $a_1=a_2+a_3+a_4+a_5$ and 
    $a_2a_3+a_2a_4+a_2a_5+a_3a_4+a_3a_5+a_4a_5=0$, 
    there exist two 1-parameter families 
    $\{x_{t,\pm}\colon{}\C\setminus\{q_1,q_2,q_3,q_4,q_5\}\to {\Real}^3
    ;\, t>0\}$ 
    of $5$-end catenoids  
    such that $\nu(q_j)=v_j$ and $w(q_j)=a_j$ $(j=1,2,3,4,5)$ .
    Indeed their representations are given by \eqref{Wf} with 
    \[
     f(z) := \frac{a_2}z
         + \frac{a_3}{z-1}
         + \frac{a_4}
                {z-\displaystyle{\frac{a_2+a_5\zeta_6}{a_2+a_3+a_5}}}
         + \frac{a_5}
                {z-\displaystyle{\frac{a_2+a_4\bar\zeta_6}{a_2+a_3+a_4}}},
    \]
    up to rigid motion in ${\Real}^3$, 
    where $\zeta_6$ is a primitive root 
    of the equation $z^6=1$,
    i.e.~$\zeta_6=(1\pm\sqrt{3}i)/2$.
    Remark here that 
    $x_{t,+}$ and $x_{t,-}$ lie on the symmetric positions 
    with respect to the $x_1x_3$-plane, 
    and hence they are isometric with each other.
\end{Exa}

We also have the following

\begin{Exa}\label{exa:4-10}
    Let $m$ be an integer greater than $1$, 
    and $\zeta_m$ a primitive root 
    of the equation $z^m=1$.
    For any positive number $t$, 
    define the surface 
    $x_t\colon{}\C\setminus\{\infty,0,1,\zeta_m,\ldots,\zeta_m^{m-1}\}
    \to {\Real}^3$  
    by the data \eqref{Wf} with 
    \[
        f(z) := \frac{(m+1)z^m+(m-1)}{z(z^m-1)}.
    \]
    Then $\{x_t\}$ is a 1-parameter family of 
    $\Z_m$-invariant $(m+2)$-end catenoids 
    with the same flux data.
    Therefore the estimate as in Theorems \ref{thm:3-3} and \ref{thm:4-2} 
    does not hold for $n$-end catenoids of TYPE I 
    for any $n\ge 4$.
\end{Exa}


\end{document}